\pdfoutput=1
\documentclass[11pt]{article}
\usepackage{bbm}
\usepackage{amsmath}
\usepackage{amsthm}
\usepackage{amssymb} 
\usepackage{mathrsfs}
\usepackage[dvipsnames]{xcolor}
\usepackage{graphicx}
\usepackage{authblk}
\usepackage[backref=none,hypertexnames=false,colorlinks=true,urlcolor=blue,linkcolor=blue,citecolor=blue]{hyperref}
\usepackage[numbers,comma,square,sort&compress]{natbib}
\usepackage[letterpaper,text={6.5in,9in},centering]{geometry}
\usepackage{tabularx}
\usepackage{booktabs} 
\usepackage{array} 
\providecommand{\keywords}[1]{\textbf{Keywords:} #1}

\graphicspath{{eps/}{pdf/}}

\makeatletter\@addtoreset{equation}{section}\makeatother

\newtheorem{thm}{Theorem}[section]

\newtheorem{rmk}{Remark}
\newtheorem{assumption}{Assumption}

\begin{document}

\title{
Bounding escape rates and approximating quasi-stationary distributions of Brownian dynamics
}
\author{Jason J. Bramburger}

\affil{\small Department of Mathematics and Statistics, Concordia University, Montr\'eal, QC, Canada}

\date{}
\maketitle

\begin{abstract}
Throughout physics Brownian dynamics are used to describe the behaviour of molecular systems. When the Brownian particle is confined to a bounded domain, a particularly important question arises around determining how long it takes the particle to encounter certain regions of the boundary from which it can escape. Termed the first passage time, it sets the natural timescale of the chemical, biological, and physical processes that are described by the stochastic differential equation. Probabilistic information about the first passage time can be studied using spectral properties of the deterministic generator of the stochastic process. In this work we introduce a framework for bounding the leading eigenvalue of the generator which determines the exponential rate of escape of the particle from the domain. The method employs sum-of-squares programming to produce nearly sharp {\color{black} numerical} upper and lower bounds on the leading eigenvalue, while also giving good {\color{black} numerical} approximations of the associated leading eigenfunction, the quasi-stationary distribution of the process. To demonstrate utility, the method is applied to prototypical low-dimensional problems from the literature.  
\end{abstract}

\keywords{stochastic differential equation, sum-of-squares, quasi-stationary distribution, Witten Laplacian, first passage time, Brownian dynamics}

\section{Introduction}

In this paper we propose numerical methods to novelly analyze the behaviour of stochastic gradient flows. Throughout we will consider a bounded domain $\Omega \subset \mathbb{R}^d$, $d \geq 1$, and a stochastic process $\{X(t)\}_{t\geq 0}$ that diffuses in $\overline\Omega$ according to the overdamped Langevin equation \cite{schlick2010molecular,ermak1978brownian}
\begin{equation}\label{SDE}
    \mathrm{d}X(t) = -\nabla V(X(t))\mathrm{d}t + \sqrt{2\sigma}\mathrm{d}B(t),
\end{equation}
where $V:\mathbb{R}^d \to \mathbb{R}$ is a smooth potential function, $B(t)$ is a $d$-dimensional Brownian motion, and $\sigma > 0$ is the noise strength. Furthermore, the boundary of the domain, denoted $\partial \Omega$, is broken into two distinguished regions: the absorbing part, $\Gamma_A$, and the reflecting part, $\Gamma_R := \partial\Omega\setminus\Gamma_A$. Conceptually, the idea is that a Brownian particle moves around in the potential landscape of $V(x)$ on the domain $\Omega$ and can only exit through the portion of the boundary attributed to $\Gamma_A$, while it is reflected back into the domain should it encounter $\Gamma_R$. Of particular importance in this problem setting is determining the {\em first passage time} \cite{redner2001guide,madrid2020competition,holcman2014time}, defined as the time it takes for a particle initiated inside $\Omega$ to leave it through $\Gamma_A$. The first passage time determines the natural timescale of many chemical, physical, and biological processes, including ions searching for membrane channels, sperm cells reaching an egg, or predators stalking their prey \cite{chou2014first}.

The stochasticity of the model equations \eqref{SDE} means that the first passage time is a random variable conditioned on the initial condition $X(0) = x$. As we will show in detail in Section~\ref{sec:FirstPassage} below, the exponential rate at which the particle escapes $\Omega$ through the boundary $\Gamma_A$ is the principal eigenvalue of the generator of the stochastic process \eqref{SDE}, given by \cite{pavliotis2014stochastic}
\begin{equation}\label{Generator}
    \mathcal{L}u := -\nabla V\cdot \nabla u  + \sigma \Delta u,   
\end{equation}
where $\Delta$ is the Laplacian operator on $\Omega$, and functions $u$ satisfy zero Dirichlet boundary conditions on $\Gamma_A$ and zero Neumann conditions on $\Gamma_R$. {\color{black} Operator $-\mathcal{L}$ is a positive elliptic operator, which in turn gives that its spectrum is real, positive, and countable with smallest eigenvalue $\lambda_0$ termed the principal eigenvalue \cite{gilbarg1977elliptic,mclean2000strongly}. If $u_0$ is the eigenfunction corresponding to the principal eigenvalue $\lambda_0$ of $-\mathcal{L}$, meaning $\mathcal{L}u_0 = - \lambda_0u_0$, standard theory of elliptic partial differential equations dictates that $\lambda_0$ is an isolated eigenvalue and $u_0$ does not vanish anywhere in $\Omega$.} The sign-definiteness of $u_0$ allows it to be interpreted (after normalization by a constant) as a probability density function, termed the {\em quasi-stationary distribution} (QSD) \cite{lelievre2024spectral,le2012mathematical,di2016jump,lelievre2022eyring}. If the initial condition $X(0)$ is drawn from the QSD then the first passage time is an exponentially distributed random variable with parameter $\lambda_0$, giving an expected first passage time of $1/\lambda_0$. Moreover, the normal derivatives of $u_0$ along $\partial \Omega$, after renormalization, give the density of the first exit points, and can thus be used to determine probabilities of exiting through each of the connected components of $\Gamma_A$.        

The QSD approach to the first passage time provides a deterministic method to analyze a stochastic problem of significant application. Thus, to approach the problem of estimating $(\lambda_0,u_0)$ this paper introduces a method for bounding $\lambda_0$ from both above and below, while also providing accurate approximations of the QSD $u_0$. Lower bounds are obtained from a variational characterization of $\lambda_0$ and follows the general method laid out in \cite{chernyavsky2023convex}, which extended previous work from \cite{valmorbida2015stability}. {\color{black} This approach is the convex dual to a previously available linear formulation using occupation measures \cite{henrion2020moment}, which has been used for nonlinear optimal control of differential equations \cite{lasserre2008nonlinear,korda2022moments}.} The method in \cite{chernyavsky2023convex} is optimally suited for the task at hand because in the specific setting of bounding leading eigenvalues of elliptic operators such as $-\mathcal{L}$, optimizers of the associated variational problem contain information about the leading eigenfunction as well. Thus, the provable existence of sequences of lower bounds converging to $\lambda_0$ lead to better approximations of $u_0$, which in turn can be used to provide upper bounds on $\lambda_0$ through the Rayleigh quotient. The result, as we will see in this manuscript, are tight lower and upper bounds on $\lambda_0$, often matching to at least 3 significant digits in our demonstrations, as well as good approximations of the QSD $u_0$.  

The presented methods can be implemented analytically to provide rigorous upper and lower bounds on $\lambda_0$. However, doing so would typically lead to conservative bounds that often come from comparing the problem on-hand with a simple one for which the solution is known, as is done, for example, in \cite{banuelos2023bounds}. The strength of the presentation that follows in this work is that for polynomial potentials $V$ and domains $\Omega$ described by polynomial inequalities, the bounding procedure can be handed over to a computer to solve a semidefinite program by appealing to the theory of polynomial optimization and sum-of-squares programming \cite{lasserre2007sum,parrilo2020sum,nesterov2000squared,parrilo2000structured,lasserre2001global}. The result is an increasing sequence of lower bounds that converge (analytically and numerically) to $\lambda_0$ as the polynomial degree of the program is increased. 

{\color{black} A major demonstration of the utility of these methods comes in the form of bounding $\lambda_0$ and approximating the QSD $u_0$ for a Brownian particle ($V \equiv 0$) attempting to escape from a circular domain in $\mathbb{R}^2$. When the exits are very small this problem is termed the {\em narrow escape problem} \cite{grebenkov2019full,holcman2004escape,schuss2007narrow,grigoriev2002kinetics,benichou2008narrow} and is well-studied in the mathematical literature using tools from asymptotic analysis \cite{cheviakov2010asymptotic,pillay2010asymptotic,delgado2015conditional,iyaniwura2021simulation,lelievre2024spectral,banuelos2023bounds}. Here we complement this analysis with sharp numerical bounds on $\lambda_0$ in the non-asymptotic regime where the common analytical tools are not applicable. Although the resulting bounds presented in this manuscript are purely numerical, there is an existing literature that employs rigorous numerics and computer-assisted proofs to bound eigenvalues of self-adjoint operators; see \cite[Chapter~10]{nakao2019numerical} and \cite{cances2020guaranteed,liu2024guaranteed}. The positioning of the bounding framework herein into a sum-of-squares problem presents a significantly different approach from these rigorous methods and, as discussed in the final section of this paper, can be extended to make the numerical method herein rigorous. Thus, the result of this paper is a novel computational bounding method that is well-founded in theory, can be adapted to a wide variety of stochastic differential equations, and numerically bound $\lambda_0$ inside a small interval.     }


This paper is organized as follows. We begin in Section~\ref{sec:FirstPassage} with a mathematical description of the problem. This includes introducing the Fokker--Planck equation and the QSD approach to the first passage time. Section~\ref{sec:Bounds} presents both the analytical and numerical bounding frameworks. In particular, \S~\ref{subsec:PDR} provides the theory for bounding $\lambda_0$ from below using its variational characterization and null Lagrangians, while \S~\ref{subsec:SOS} describes how to translate this problem into a sum-of-squares program to be solved numerically. \S~\ref{subsec:UpperBnd} then shows how information from the numerical lower bounds can be leveraged to produce numerical upper bounds using the Rayleigh quotient. Section~\ref{sec:Demonstrations} provides numerical demonstrations, including characterizing the asymptotic escape rate from a double well potential and the aforementioned application to a Brownian particle escaping a circular domain. We conclude in Section~\ref{sec:Discussion} with a discussion of the results, avenues for future work, and directions towards validating these numerical results to make them rigorous.

\section{First Passage Times}\label{sec:FirstPassage}

Recall from the introduction that we are interested in the stochastic differential equation \eqref{SDE} restricted to a bounded and open set $\Omega \subset \mathbb{R}^d$, $d \geq 1$. Moreover, the boundary is decomposed into disjoint absorbing and reflecting components $\partial\Omega = \Gamma_A\sqcup \Gamma_R$ with $\Gamma_A \neq \emptyset$. {\color{black} We further assume that $\Gamma_A$ is comprised of finitely-many smooth, disjoint, and nonempty components sufficient to define a trace operator on; see \cite{mclean2000strongly}.} Then, for an initial condition $X(0) = x \in \Omega$ and $\{X(t)\}_{t \geq 0}$ satisfying \eqref{SDE}, we define the first passage time $\tau > 0$ as 
\begin{equation}
    \tau := \inf\{t \geq 0:\ X(t)\in \Gamma_A\}.
\end{equation}
Notice that $\tau$ is a random variable that depends on the initial condition, so we can define the {\em survival probability} conditioned on the starting position as
\begin{equation}
    S(x,t) := \mathbb{P}(\tau > t|\ X(0) = x).
\end{equation}
The survival probability satisfies the Fokker--Planck (or Kolmogorov backward) equation \cite{madrid2020competition,pavliotis2014stochastic}
\begin{equation}\label{SurvivalPDE}
    \begin{split}
        \frac{\partial S}{\partial t} &= \mathcal{L} S, \quad x \in \Omega, \\
        S &= 0, \quad x \in \Gamma_A, \\
        \vec{n}\cdot \nabla S &= 0, \quad x\in\Gamma_R,
    \end{split}
\end{equation}
where $\mathcal{L}$ is the generator \eqref{Generator} and $\vec{n}:\partial\Omega \to \mathbb{R}^d$ denotes the unit outward normal to $\Omega$. The evolution of \eqref{SurvivalPDE} is initialized with $S(x,0) = 1$ since if $x\in\Omega$, then $\tau > 0$ with probability 1.  

We introduce the Boltzmann-type weight function 
\begin{equation}\label{WeightFn}
    w(x) = c\mathrm{e}^{-V(x)/\sigma},
\end{equation}
where $c > 0$ is a normalization constant. {\color{black} The operator $\mathcal{L}$ is non-positive and symmetric on the weighted space 
\begin{equation}\label{WeightedL2}
    L^2_w(\Omega) = \{f:\Omega \to \mathbb{R}:\ \int_\Omega |f(x)|^2w(x)\mathrm{d}x < \infty\}, 
\end{equation}
equipped with the inner product 
\begin{equation}
    \langle f,g\rangle_w := \int_\Omega f(x)g(x)w(x)\mathrm{d}x.
\end{equation}
Furthermore, the non-negative operator $-\mathcal{L}$ can be defined as the Friedrichs extension of the quadratic form $Q(f) = \int_\Omega |\nabla f(x)|^2 w(x)\mathrm{d}x$ with form domain 
\begin{equation}
    \mathcal{D}(Q) = H^1_{0,\Gamma_A}(\Omega) := \{f \in H^1(\Omega):\ f|_{\Gamma_A} = 0\}, 
\end{equation}
the space of $H^1(\Omega)$ functions whose trace vanishes on $\Gamma_A$. The domain of $\mathcal{L}$ is not $H^2(\Omega)$ due to the mixed boundary conditions (see \cite[Remark~2.2]{lelievre2024spectral}), but the subspace 
\begin{equation}
    \mathcal{D}(\mathcal{L}) = \{f \in H^1_{0,\Gamma_A}(\Omega):\ \Delta f \in L^2(\Omega),\ f|_{\Gamma_A} = 0,\ \vec{n}\cdot \nabla f|_{\Gamma_R} = 0\},
\end{equation}
where $L^2(\Omega) := L^2_1(\Omega)$ in the notation of \eqref{WeightedL2}, i.e.  a weight function of $w \equiv 1$. As an aside, since $\Omega$ is bounded and $V$ is smooth, a function in $L^2_w(\Omega)$ is also in $L^2(\Omega)$ and the two norms are equivalent. The interest in introducing the innner product $\langle\cdot,\cdot\rangle_w$ is that the operator $\mathcal{L}$ is symmetric with respect to this inner product. 

The associated eigenfunctions $\{u_n\}_{n = 0}^\infty$ of $\mathcal{L}$ satisfy
\begin{equation}\label{SurvivalEigenfunctions}
    \begin{split}
        -\mathcal{L}u_n &= \lambda_nu_n, \quad x \in \Omega, \\
        u_n &= 0, \quad x \in \Gamma_A, \\
        \vec{n}\cdot \nabla u_n &= 0, \quad x\in\Gamma_R,
    \end{split}
\end{equation}
with a countably infinite sequence of positive eigenvalues \cite[Theorem~4.12]{mclean2000strongly}
\begin{equation}\label{EigenSequence}
    0 < \lambda_0 < \lambda_1 \leq \lambda_2 \leq \cdots \leq \lambda_n \leq \lambda_{n+1} \leq \cdots.
\end{equation}
Notice that the eigenvalue-eigenfunction pairs in \eqref{SurvivalEigenfunctions} can be interpreted in two ways: (i) $(\lambda_n,u_n)$ are eigenpairs of the non-negative operator $-\mathcal{L}$, or (ii) $(-\lambda_n,u_n)$ are eigenpairs of the non-positive operator $\mathcal{L}$. Both formulations are equivalent, while the election to present the eigenpairs in this way is to keep the eigenvalues positive in order to easily see the correspondence with the QSD and the first passage time $1/\lambda_0$.   
Regardless of the interpretation, the eigenfunctions form an orthonormal basis for $L^2_w(\Omega)$ and so we can expand $S(x,t)$ as a solution of \eqref{SurvivalPDE} in the eigenfunctions to obtain }
\begin{equation}\label{EigenfunctionExpansion}
    S(x,t) = \sum_{n = 0}^\infty s_n u_n(x)\mathrm{e}^{-\lambda_n t},
\end{equation}
where $s_n = \langle u_n,1\rangle_w$. Supposing further that the initial condition $x$ is chosen from a probability measure $\mu$, then the survival probability conditioned on the initial distribution is given by
\begin{equation}
    \mathbb{P}(\tau > t|\ X(0) \sim \mu) = \int_\Omega S(x,t)\mathrm{d}\mu(x).
\end{equation}
Owing to the strict inequalities $0 < \lambda_0 < \lambda_1$, we have the exponential rate
\begin{equation}
    \lambda_0 = -\lim_{t \to \infty} \frac{1}{t}\log\mathbb{P}(\tau > t|\ X(0) \sim \mu)
\end{equation}
at which the particle hits the boundary $\Gamma_A$. 

Standard elliptic partial differential equation theory gives that the leading eigenfunction $u_0$ is nonzero over all of $\Omega$. The QSD for the process is given by
\begin{equation}
    \mathrm{d}\mu_0(x) = \frac{u_0(x) w(x) \mathrm{d} x}{\int_\Omega u_0(y) w(y)\mathrm{d}y},
\end{equation}
where $w$ is the weight function \eqref{WeightFn}. If the initial position of the particle is drawn from the probability measure $\mu_0$ we have 
\begin{equation}
    \mathbb{P}(\tau > t|\ X(0) \sim \mu_0) = \mathrm{e}^{-\lambda_0 t}  
\end{equation}
for all $t \geq 0$, which can be easily checked from \eqref{EigenfunctionExpansion}. Hence, $\tau$ is now an exponential random variable with $\mathbb{E}[\tau] = \lambda_0^{-1}$. While it might not be expected that one can draw from the QSD directly, \eqref{EigenfunctionExpansion} gives that $S(x,t) \approx s_0 u_0(x)\mathrm{e}^{-\lambda_0 t}$ for $t \gg 1$. Thus, if the stochastic process stays trapped for a long time before leaving, such as when $\Gamma_A$ is small relative to $\Gamma_R$, then starting from the QSD approximates the long-time behaviour of the stochastic process reaching the QSD before leaving the domain $\Omega$. Furthermore, if $\Gamma_A$ is the disjoint union of $k\geq 1$ connected components $A_1,\dots,A_k \subset \partial\Omega$, then we have the probabilities (see for example the distribution $\rho$ in the proof of \cite[Proposition 2.4]{le2012mathematical})
\begin{equation}
    \mathbb{P}(X(\tau) \in A_j|\ X(0) \sim \mu_0) = \frac{\int_{A_j} \vec{n}\cdot \nabla u_0(x) w(x) \mathrm{d} x}{\int_{\Gamma_A} \vec{n}\cdot \nabla u_0(x) w(x) \mathrm{d} x}, \qquad j = 1,\dots, k, 
\end{equation}
describing the probability of exiting $\Omega$ through each of the components $A_j$. Finally, one can check that $\tau$ and $X(\tau)$ are independent random variables, and thus that they can be sampled using a jump Markov model; see \cite{di2016jump} for details.


\section{Bounds on the Principal Eigenvalue}\label{sec:Bounds}

In this section we provide the analytical and computational bounding framework for {\color{black}the principal eigenvalue $\lambda_0$, as defined in \eqref{EigenSequence}}, with associated leading eigenfunction $u_0$. We first symmetrize the elliptic operator $\mathcal{L}$ by introducing the unitary transformation of $\mathcal{L}$ to
\begin{equation}\label{unitarytransformation}
    \mathcal{H} = w^{1/2}\mathcal{L}w^{-1/2}. 
\end{equation}
The result is that we have re-cast $\mathcal{L}$ as a time-independent Schr\"odinger equation called the Witten Laplacian \cite{michel2019small,witten1982supersymmetry}, given by 
\begin{equation}
    \mathcal{H} = \sigma\Delta - U(x),
\end{equation}
where 
\begin{equation}
    U(x) = \frac{|\nabla V(x)|^2}{\sigma} - \frac{\Delta V}{2}
\end{equation}
and $V$ is the potential from \eqref{SDE}. The domain of the operator $\mathcal{H}$ is defined as the image of the domain of the operator $\mathcal{L}$ under the unitary transformation, i.e. $\mathcal{D}(\mathcal{H}) = w^{1/2}\mathcal{D}(\mathcal{L})$. Then $(\lambda,u)$ is an eigenvalue-eigenfunction pair for $\mathcal{L}$ if and only if $(\lambda,uw^{1/2})$ is an eigenvalue-eigenvector pair for $\mathcal{H}$. However, the functions in $ \mathcal{D}(\mathcal{H})$ do not necessarily obey the same boundary conditions as those in $\mathcal{D}(\mathcal{L})$. Indeed, if $v = w^{1/2}u \in \mathcal{D}(\mathcal{H})$ with $u \in \mathcal{D}(\mathcal{L})$, then $v$ vanishes on $\Gamma_A$ if and only if $u$ vanishes on $\Gamma_A$ as well since $w > 0$ everywhere, but the Neumann conditions imposed on $\Gamma_R$ for $u$ become Robin boundary conditions for $v$ since we get
\begin{equation}
          \vec{n}\cdot \nabla u = \vec{n}\cdot \nabla (w^{-1/2}v) = 0 \implies \vec{n}\cdot\nabla v + \frac{v}{2\sigma}\vec{n}\cdot\nabla V  = 0, \qquad \forall x \in \Gamma_R,
\end{equation}
again using that $w > 0$. Since Robin boundary conditions do not fit within the computational framework outlined in this section, we make the following simplifying assumption to proceed through this section.

\begin{assumption}\label{assump1}
    The potential $V:\mathbb{R}^d \to \mathbb{R}$ is such that for all $x \in \Gamma_R$ we have $\vec{n}\cdot \nabla V = 0$.
\end{assumption}

Although this assumption is restrictive, it can be satisfied in at least three relevant scenarios:
\begin{enumerate}
    \item fully absorbing boundaries, meaning $\Gamma_R = \emptyset$;
    \item Brownian motion in a homogeneous landscape, meaning $V \equiv 0$;
    \item the domain $\Omega$ is taken to be the boundary of a basin of attraction of a local minimum of a potential $V$ for the deterministic gradient flow dynamics $\dot x = -\nabla V(x)$.
\end{enumerate}
Indeed, (1) is the setting for our convergence results that follow, while demonstrations coming from notable problems in the literature for cases (1) and (2) above are presented in Section~\ref{sec:Demonstrations}. The result of Assumption~\ref{assump1} is that the functions in $\mathcal{D}(\mathcal{H})$ now satisfy the same mixed Dirichlet--Neumann boundary conditions as those in $\mathcal{D}(\mathcal{L})$. 

To bound $\lambda_0 > 0$ from below we will use the variational characterization of the principal eigenvalue of the non-negative operator $-\mathcal{H}$:
\begin{equation}\label{RQ}
    \lambda_0 = \sup_{\lambda \in \mathbb{R}}\ \lambda\quad \mathrm{s.t.}\quad \inf_{u \in H^1_{0,\Gamma_A}(\Omega)}\int_\Omega [\sigma|\nabla u|^2 + (U(x) - \lambda )u^2] \mathrm{d}x \geq 0,
\end{equation}
where the reflecting boundary condition $\vec{n}\cdot \nabla u|_{\Gamma_R} = 0$ inherited by Assumption~\ref{assump1} is naturally included in the variational formulation through the Euler--Lagrange equation associated to this optimization problem. The variational formula \eqref{RQ} for $\lambda_0$ follows from a rearrangement of the Rayleigh quotient associated to the elliptic operator $-\mathcal{H}$. In \S~\ref{subsec:PDR} we introduce the method from \cite{chernyavsky2023convex} that makes use of \eqref{RQ} to bound $\lambda_0$ from below. Then, in \S~\ref{subsec:SOS} we translate the bounding problem into one of polynomial optimization which can be implemented computationally as a sum-of-squares program. Finally, in \S~\ref{subsec:UpperBnd} we show how information from the sum-of-squares lower bound can be used to bound $\lambda_0$ from above using the Rayleigh quotient and approximate the QSD.

\subsection{Lower Bounds on the Principal Eigenvalue}\label{subsec:PDR}

In this subsection we describe a method of bounding the principal eigenvalue, $\lambda_0$, of the operator $-\mathcal{L}$ from below using the variational formulation \eqref{RQ}. To begin, let $\varphi:\Omega \times \mathbb{R} \to \mathbb{R}^d$ be a function for which 
\begin{equation}
    \int_\Omega \nabla \cdot \varphi(x,u) \mathrm{d}x =  \int_\Omega \bigg[\sum_{i = 1}^d \frac{\partial \varphi}{\partial x_i} + \frac{\partial \varphi}{\partial u}(\nabla\cdot u)\bigg] \mathrm{d}x = 0,
\end{equation}
{\color{black}for all $u \in H^1_{0,\Gamma_A}(\Omega)$.} Any $\varphi(x,u)$ satisfying the above is a null Lagrangian and trivially gives that 
\begin{equation}
    \int_\Omega [\sigma|\nabla u|^2 + (U(x) - \lambda )u^2]  \mathrm{d}x = \int_\Omega [\sigma|\nabla u|^2 + (U(x) - \lambda )u^2 + \nabla \cdot \varphi(x,u)] \mathrm{d}x.
\end{equation}
We estimate the above integral from below by taking the pointwise infimum of its integrand over all $(x,u,\nabla u)$. This leads to the inequality
\begin{equation}\label{PointwiseRelax}
    \begin{split}
        \int_\Omega [\sigma|\nabla u|^2 + &(U(x) - \lambda )u^2 + \nabla \cdot \varphi(x,u)] \mathrm{d}x \\ 
        &\geq \int_\Omega \inf_{\substack{y \in \mathbb{R} \\ z \in \mathbb{R}^d}} \bigg[\sigma|z|^2 +  (U(x) - \lambda )y^2 + \sum_{i = 1}^d \bigg(\frac{\partial \varphi}{\partial x_i}(x,y) + \frac{\partial \varphi}{\partial u}(x,y)z_i\bigg)\bigg]\mathrm{d}x.
    \end{split}
\end{equation}
Thus, for any fixed $\varphi$ and $\lambda \in \mathbb{R}$, if 
\begin{equation}
    \sigma|z|^2 +  (U(x) - \lambda )y^2 + \sum_{i = 1}^d \bigg(\frac{\partial \varphi}{\partial x_i}(x,y) + \frac{\partial \varphi}{\partial u}(x,y)z_i\bigg) \geq 0
\end{equation}
for all $(x,y,z) \in \Omega\times \mathbb{R}\times\mathbb{R}^d$, it follows from the inequality \eqref{PointwiseRelax} that $\lambda_0 \geq \lambda$, thus providing a lower bound on the principal eigenvalue of $\mathcal{L}$.

The above process can be taken further by optimizing over functions $\varphi$. The result is the {\em pointwise dual relaxation} method \cite{chernyavsky2023convex} which returns a lower bound $\lambda_\mathrm{pdr} \leq \lambda_0$ via the optimization problem
\begin{equation}\label{PDR}
   \lambda_\mathrm{pdr} := \sup_{\varphi, \lambda}\ \lambda\quad \mathrm{s.t.}\quad \sigma|z|^2 +  (U(x) - \lambda )y^2 + \sum_{i = 1}^d \bigg(\frac{\partial \varphi}{\partial x_i}(x,y) + \frac{\partial \varphi}{\partial u}(x,y)z_i\bigg) \geq 0,\ \forall (x,y,z) \in \Omega\times\mathbb{R}\times\mathbb{R}^d.
\end{equation}
The following theorem demonstrates that there is no relaxation gap, i.e. $\lambda_\mathrm{pdr} = \lambda_0$, in the case of purely Dirichlet boundary conditions ($\Gamma_R = \emptyset$).

\begin{thm}\cite[Theorem~4.1(1)]{chernyavsky2023convex}\label{thm:PDRConv}
    Let $\Omega$ be an open, bounded, and Lipschitz domain, $U:\mathbb{R}^d \to \mathbb{R}$ a smooth function, $\sigma > 0$, and $\Gamma_A = \partial \Omega$. Then, $\lambda_0 = \lambda_\mathrm{pdr}$. In particular, there exists $\lambda^k \in \mathbb{R}$ and $f^k\in C(\overline\Omega,\mathbb{R}^d)$ such that, with $\varphi^k = f^k(x)u^2$, the sequence $\{\lambda^k,\varphi^k\}_{k = 1}^\infty$ satisfies
    \begin{equation}\label{PDRInequality}
        \sigma|z|^2 +  (U(x) - \lambda^k )y^2 + \sum_{i = 1}^d \bigg(\frac{\partial \varphi^k}{\partial x_i}(x,y) + \frac{\partial \varphi^k}{\partial u}(x,y)z_i\bigg) \geq 0, \quad \forall (x,y,z) \in \Omega\times\mathbb{R}\times\mathbb{R}^d    
    \end{equation}
    for all $k \geq 1$ and $\lambda^k \nearrow \lambda_0$ as $k \to \infty$.
\end{thm}

Notice that with $\varphi(x,u) = f(x)u^2$ as in Theorem~\ref{thm:PDRConv}, the divergence theorem gives that
\begin{equation}
    \int_\Omega \nabla\cdot f(x)u^2\mathrm{d}x = \int_{\partial\Omega} \vec{n}\cdot f(x)u^2 \mathrm{d}S. 
\end{equation}
In the case of Dirichlet boundary conditions, i.e. $\Gamma_A = \partial\Omega$, $u$ vanishes everywhere on the boundary and so $\varphi(x,u) = f(x)u^2$ is indeed a null Lagrangian for any continuous function $f$. {\color{black} While the full proof of Theorem~\ref{thm:PDRConv} can be found in \cite{chernyavsky2023convex}, we provide a brief sketch since it will be used to generate the upper bounds on $\lambda_0$ in Section~\ref{subsec:UpperBnd} below.

\begin{proof}[Sketch of Proof of Theorem~\ref{thm:PDRConv}]
    Putting the null Lagrangian $\varphi(x,u) = f(x)u^2$ into the desired inequality \eqref{PDRInequality} gives
\begin{equation}
        \sigma|z|^2 + (U(x) - \lambda )y^2 + (\nabla\cdot f(x)) y^2 + 2yf(x)\cdot z \geq 0. 
\end{equation}
The above inequality can be equivalently rewritten as the bilinear form
\begin{equation}
        \begin{bmatrix}
            y \\ z
        \end{bmatrix}^\intercal \begin{bmatrix}
            \nabla\cdot f(x) + (U(x) - \lambda ) & f(x)^\intercal \\ f(x) & \sigma I_d   
        \end{bmatrix}\begin{bmatrix}
            y \\ z
        \end{bmatrix} \geq 0,
\end{equation}
where $I_d$ is the $d\times d$ dimensional identity matrix. The bilinear form being non-negative is equivalent to having 
\begin{equation}
    \begin{bmatrix}
            \nabla\cdot f(x) + (U(x) - \lambda ) & f(x)^\intercal \\ f(x) & \sigma I_d    
        \end{bmatrix} \succeq 0, \quad \forall x \in \Omega,
\end{equation}
where $\succeq 0$ denotes the matrix being positive semidefinite. Using the Schur complement, this now becomes equivalent to
\begin{equation}\label{PDRu1}
    \nabla\cdot f(x) + (U(x) - \lambda ) - \sigma^{-1} f(x)\cdot f(x) \geq 0.   
\end{equation}   
Now, the specific choice $f(x) = -\sigma\nabla u_0/u_0$, where $u_0$ is the principal eigenfunction of $-\mathcal{H}$ associated to $\lambda_0$, will reduce \eqref{PDRu1} to $\lambda_0 - \lambda \geq 0$. Maximizing $\lambda_0 - \lambda \geq 0$ over $\lambda$ thus results in the principal eigenvalue $\lambda_0$.

However, one should note that the Dirichlet boundary conditions give that $f(x) = -\sigma \nabla u_0/u_0$ is not a valid choice for the function $f(x)$ since it is not bounded (and not continuous) on the boundary of $\Omega$. Thus, the full proof of Theorem~\ref{thm:PDRConv} provides a sequence that approximates this optimal choice while remaining continuous and bounded up to the boundary. Roughly, this is attained by creating a sequence of larger domains $\{\Omega_k\}_{k=1}^\infty$ which shrink down on $\Omega$ as $k \to \infty$ and obtaining principal eigenfunctions $u_0^k$ to $-\mathcal{H}$ on the $\Omega_k$ domains with Dirichlet boundary conditions. The proof then takes $f^k(x) = -\sigma \nabla u_0^k(x)/u_0^k(x)$ and exploits the continuity of the spectrum of $-\mathcal{H}$ with respect to the domain to obtain the convergence results. This completes the sketch. 
\end{proof}

\begin{rmk}\label{rmk:GammaR}
    Theorem~\ref{thm:PDRConv} is specific to the case of Dirichlet boundary conditions. However, the above sketch can be generalized to the case of mixed boundary conditions, i.e. $\Gamma_A,\Gamma_R \neq \emptyset$, by again taking $\varphi(x,u) = f(x)u^2$ with $f(x) = -\sigma\nabla u_0(x)/u_0(x)$. One finds that this choice is again a null Lagrangian thanks to the fact that $\vec{n}\cdot f(x) = 0$ for all $x \in \Gamma_R$ and that \eqref{PDRu1} will again reduce to $\lambda_0 - \lambda \geq 0$. The technical hurdle in making this sketch rigorous comes from proving the continuity of the spectrum of $-\mathcal{H}$ with mixed boundary conditions with respect to the domain to complete the proof as outlined above for the Dirichlet case. The author was unable to find general results on this domain continuity and therefore the theoretical results are left specific to the purely Dirichlet case, while our numerical results below indicate that these proofs (in certain cases) could likely be achieved.    
\end{rmk}

}

\subsection{Relaxation to a Sum-of-Squares Program}\label{subsec:SOS}

In this subsection we turn to a computational implementation for extracting the convergent sequences guaranteed by Theorem~\ref{thm:PDRConv}. In particular, we will demonstrate how to incorporate sum-of-squares (SOS) polynomial optimization into the problem, which can be reformulated as a semidefinite program (SDP) to be solved numerically using interior point methods. 

We begin by assuming that $U(x)$ in \eqref{RQ} is a polynomial in $x$ and the set $\overline\Omega$ is a basic closed semialgebraic set, meaning that there exists polynomials $g_1(x),\dots g_s(x)$ so that 
\begin{equation}\label{OmegaSemiAlg}
    \overline\Omega = \{x\in\mathbb{R}^d:\ g_1(x) \geq 0, \cdots, g_s(x) \geq 0\}.
\end{equation}
Suppose further that the reflecting boundary component $\Omega$ can be broken into $r$ disjoint components $\Gamma_R = \Gamma_1 \sqcup \cdots \sqcup \Gamma_r$ defined by a single $g_i(x)$ vanishing on component $\Gamma_i$, i.e.
\begin{equation}
    \Gamma_i = \{x\in\mathbb{R}^d:\ g_i(x) = 0\ \mathrm{and}\ g_j(x) \geq 0, j \neq i\},
\end{equation}
and that the gradients $\nabla g_i(x)$ do not vanish anywhere on $\Gamma_R$ for all $i = 1,\dots, r$. These assumptions will allow for the formulation of \eqref{PDRInequality} as a polynomial optimization problem that can be strengthened to an SOS condition to obtain numerically converging lower bounds on $\lambda_0$. We may further generalize the above, including to $\Omega$ being the disjoint union of semialgebraic sets, as is demonstrated in \S~\ref{subsec:NarrowEscape} below, but refrain from doing so here to simplify the presentation. 

Let $\mathbb{R}^d[x]_\nu$ be the set of $d$-dimensional vectors of polynomials in $x$ of maximal degree $\nu$. Then, taking $\varphi(x,u) = f(x)u^2$ with $f \in \mathbb{R}^d[x]_\nu$ gives the trivial lower bound 
\begin{equation}\label{PDRpolyopt}
   \lambda_\mathrm{pdr} \geq \lambda_{\nu} := \sup_{\substack{\lambda\in\mathbb{R}\\f\in\mathbb{R}^d[x]_\nu}}\ \lambda\quad \mathrm{s.t.}\quad \begin{cases}
        \sigma|z|^2 + (U(x) - \lambda )y^2 + (\nabla\cdot f(x)) y^2 + 2yf(x)\cdot z \geq 0\ \mathrm{on}\ \Omega\times\mathbb{R}\times\mathbb{R}^d, \\
        \nabla g_i(x)\cdot f(x) = 0\ \mathrm{on}\ \Gamma_i, \ i = 1,\dots,r.
   \end{cases}
\end{equation}
The optimization problem \eqref{PDRpolyopt} is finite-dimensional for any $\nu \in \mathbb{N}$ and generates an increasing sequence of lower bounds by increasing $\nu$. {\color{black}Unfortunately, verifying pointwise inequalities of multi-dimensional polynomials is an NP-hard problem in general \cite{murty1987some}. We therefore proceed using a standard approach to strengthen all polynomial inequalities in \eqref{PDRpolyopt} to be weighted SOS constraints \cite{lasserre2001global}. The result will be a convex maximization problem in which the polynomial coefficients appear only linearly in the constraints of the problem which can be reformulated as an SDP to be solved numerically. } 

Let us denote $\Sigma[\cdot]$ to be the set of polynomials in the bracketed variables that can be written as sums of squares of other polynomials. We consider the sets 
\begin{equation}\label{QOmega}
    \begin{split}
        \mathcal{Q}_\nu(\Omega) &= \bigg\{\rho_0 + \sum_{i = 1}^s g_i\rho_i:\ \rho_k \in \Sigma[x,y,z]\cap \mathbb{R}[x,y,z]_\nu\bigg\} \\
        \mathcal{Q}_\nu(\Gamma_i) &= \bigg\{\zeta g_i + \rho_0 + \sum_{j \neq i}^s g_j\rho_j:\ \rho_k \in \Sigma[x]\cap \mathbb{R}[x]_\nu,\ \zeta \in\mathbb{R}[x]_\nu\bigg\}. 
    \end{split}
\end{equation}
Polynomials that belong to $\mathcal{Q}_\nu(\Omega)$ are non-negative on $\Omega$ since \eqref{OmegaSemiAlg} gives that $g_i \geq 0$ on $\Omega$ and $\rho_k \geq 0$ are SOS. Similarly, polynomials in $\mathcal{Q}_\nu(\Gamma_i)$ are non-negative on $\Gamma_i$ since $g_i = 0$, while $g_j \geq 0$, $i \neq j$ and $\rho_k \geq 0$. With these sets we may strengthen the polynomial inequalities in \eqref{PDRpolyopt} to weighted SOS constraints 
\begin{equation}\label{PDRsos}
   \lambda_{\nu} \geq \lambda_{\nu}^\mathrm{SOS} := \sup_{\substack{\lambda\in\mathbb{R}\\f\in\mathbb{R}^d[x]_\nu}}\ \lambda\quad \mathrm{s.t.}\quad \begin{cases}
        \sigma|z|^2 + (U(x) - \lambda )y^2 + (\nabla\cdot f(x)) y^2 + 2yf(x)\cdot z \in \mathcal{Q}_\nu(\Omega), \\
        \nabla g_i(x)\cdot f(x) \in \mathcal{Q}_\nu(\Gamma_i), \ i = 1,\dots,r, \\
        -\nabla g_i(x)\cdot f(x) \in \mathcal{Q}_\nu(\Gamma_i),\ i = 1,\dots,r.
   \end{cases}
\end{equation}
Notice that condition $\nabla g_i(x)\cdot f(x) \in \mathcal{Q}_\nu(\Gamma_i)$ gives that $\nabla g_i(x)\cdot f(x) \geq 0$ on $\Gamma_i$ and $-\nabla g_i(x)\cdot f(x) \in \mathcal{Q}_\nu(\Gamma_i)$ gives that $\nabla g_i(x)\cdot f(x) \leq 0$ on $\Gamma_i$, which together give the desired equality condition $\nabla g_i(x)\cdot f(x) = 0$ on $\Gamma_i$. The advantage now is that the SOS constraints in \eqref{PDRsos} can be reformulated as SDPs and solved numerically \cite{nesterov2000squared,parrilo2000structured,lasserre2001global}, thus providing the computational lower bounds on $\lambda_0 \geq \lambda_\nu^\mathrm{SOS}$.

The following theorem shows that the SOS lower bounds converge to $\lambda_0$ as $\nu \to \infty$ when $\Gamma_R = \emptyset$. 

\begin{thm}\cite[Theorem~4.1(2)]{chernyavsky2023convex}\label{thm:PDRsos}
    If $U$ is polynomial, $\Gamma_A = \partial\Omega$, and $\Omega$ is a basic semialgebraic set with $r^2 - \|x\|^2 \in \mathcal{Q}_\nu(\Omega)$ for some $\nu\in\mathbb{N}$ and $r \in \mathbb{R}$, then $\lambda_\nu^\mathrm{SOS}\nearrow\lambda_0$ as $\nu \to \infty$.
\end{thm}

Notice that the additional constraint $r^2 - \|x\|^2 \in \mathcal{Q}_\nu(\Omega)$ for some $r \in \mathbb{R}$ implies that $\Omega$ is bounded since $\|x\|^2 \leq r^2$ for all $x \in \Omega$. Unfortunately, boundedness of $\Omega$ is not always enough to guarantee this condition \cite[Chapter~2]{lasserre2015introduction}, but it can be satisfied by adding the additional constraint $g_{s+1}(x) = r_0^2 - \|x\|^2$  into \eqref{OmegaSemiAlg} with $r_0 > 0$ sufficiently large to not change the definition of the set itself. {\color{black}Finally, like Theorem~\ref{thm:PDRConv}, Theorem~\ref{thm:PDRsos} is only proven for the case that $\Gamma_R = \emptyset$. If one is able to overcome the technical hurdles outlined in Remark~\ref{rmk:GammaR} to obtain a version of Theorem~\ref{thm:PDRConv} with $\Gamma_R = \emptyset$, then one can readily apply \cite[Theorem~3.5]{chernyavsky2023convex} to prove Theorem~\ref{thm:PDRsos} with $\Gamma_R \neq \emptyset$.} In the numerical demonstrations that follow we see no measurable difference in convergence of the lower bounds to $\lambda_0$ between the $\Gamma_R = \emptyset$ and $\Gamma_R \neq \emptyset$ cases so long as the other assumptions are satisfied.

\subsection{Upper Bounds on the Principal Eigenvalue}\label{subsec:UpperBnd}

Upper bounds on $\lambda_0$ can always be attained through the Rayleigh quotient. That is, $\lambda_0 \leq R(u)$ where
\begin{equation}\label{RQ2}
    R(u) = \frac{\int_\Omega \sigma|\nabla u|^2 + U(x)u^2 \mathrm{d}x}{\int_\Omega u^2 \mathrm{d}x}
\end{equation}
{\color{black}and $u \in H^1_{0,\Gamma_A}$.} Note that the leading eigenfunction $u_0$ of $-\mathcal{H}$ will give $\lambda_0 = R(u_0)$, while any other function gives $\lambda_0 < R(u)$. Thus, tight upper bounds are achieved by good approximations of the leading eigenfunction, for which we can appeal to the theory of the previous subsections to obtain.

To begin, recall from Theorem~\ref{thm:PDRConv} and the discussion following it that near-optimizers of the pointwise dual relaxation method \eqref{PDR} take the form $\varphi(x,u) = f(x)u^2$ with 
\begin{equation}
    f(x) \approx -\frac{\sigma\nabla u_0(x)}{u_0(x)} \implies f(x) \approx -\sigma\nabla \ln(u_0(x)).
\end{equation}
Therefore, we may use the gradient theorem to approximate
\begin{equation}
    \sigma\ln(u_0(x)) \approx u_0(x_0) + \int_0^1 f(\gamma(t))\cdot\gamma'(t)\mathrm{d}t,
\end{equation}
where $\gamma(t) = x_0 + t(x - x_0)$ is a parametric curve in $\Omega$ with starting point $x_0 \in \Omega$ fixed and $u_0(x_0)$ playing the role of a constant of integration. Dividing by $\sigma$ and taking an exponential in turn approximates $u_0(x)$. 

With this approximation of the leading eigenfunction of $\mathcal{H}$, we may then also compute its gradient for use of \eqref{RQ2} to bound $\lambda_0$ above. Indeed, using the identity $u_0(x) = \mathrm{e}^{\ln(u_0(x))}$, we obtain the gradient 
\begin{equation}
    \nabla u_0(x) = \mathrm{e}^{\ln(u_0(x))}\nabla \ln(u_0(x)) \approx \frac{u_0(x)f(x)}{\sigma}.
\end{equation}
Thus, no numerical derivatives are required to approximate the gradient of the principal eigenvalue. As we will see in the examples below, this method can provide tight upper bounds on the leading eigenvalue that can match the lower bounds from the previous subsection to 3 or more significant digits. 

Finally, notice that our approximation of the leading eigenfunction for $\mathcal{H}$ can in turn lead to an approximation for the leading eigenfunction for $\mathcal{L}$. Indeed, if $u_0(x)$ is an approximation of the leading eigenfunction of $\mathcal{H}$, then the weight function $w(x)$ in \eqref{WeightFn} gives that $u_0(x)(w(x))^{1/2}$ is an approximation of the leading eigenfunction of $\mathcal{L}$.

\section{Numerical Demonstrations}\label{sec:Demonstrations}

In this section we provide numerical demonstrations of the theory from the previous section using semidefinite programming. In particular, we use YALMIP \cite{lofberg2004yalmip} to reformulate SOS problems as SDPs and MOSEK to solve them \cite{aps2019mosek}. The combination of YALMIP and MOSEK can be used to tune linear combinations of coefficients for a function to admit a representation as an SOS polynomial. In this way, we verify a polynomial $p$ belongs to $\mathcal{Q}_\nu(\Omega)$ from \eqref{QOmega} by numerically identifying monomial coefficients up to degree $\nu$ so that each polynomial in the set   
\begin{equation}
    \bigg\{p - \sum_{i = 1}^s g_i\rho_i, \rho_1,\dots,\rho_s\bigg\}
\end{equation}
admits an SOS representation. This strategy can be observed in the associated code that accompanies this manuscript, freely available \cite{gitrepo}, for which all demonstrations that follow can be reproduced in MATLAB.

\subsection{Escape From a Double Well Potential}

Consider the 1-dimensional stochastic double well problem
\begin{equation}
    \mathrm{d}X(t) = X(t)(1 - X(t))(X(t) + \alpha)\mathrm{d}t + \sqrt{2}\mathrm{d}B(t)    
\end{equation}
where $\alpha \in [0,1]$ is a parameter that moves the position of one of the potential wells. We further consider two absorbing boundaries at $x = \pm L$. In this case the generator is given by
\begin{equation}\label{DoubleWellOperator}
    \mathcal{L}u = u'' + x(1 - x)(x + \alpha)u',
\end{equation}
acting on functions satisfying the boundary conditions $u(\pm L) = 0$. The generator is symmetric in the $L^2$ space weighted by $w(x) = \mathrm{e}^{-x^2[\alpha(4x - 6) + x (3x - 4)]/(12\sigma)}$, while the symmetrized elliptic operator takes the form
\begin{equation}
    \mathcal{H}u = u'' - \underbrace{\frac{1}{2}[2x^2(1 - x)^2(x + \alpha)^2 - \alpha(2x - 1) + (2 - 3x)x]}_{= U(x)}u,
\end{equation}
again acting on functions with $u(\pm L) = 0$. To improve numerical conditioning we will rescale space as $x \mapsto Lx$, resulting in the scaled elliptic operator
\begin{equation}
    \mathcal{H}_L = \frac{1}{L^2}u'' - U(Lx)u 
\end{equation}
now acting on functions with $u(\pm 1) = 0$ and $L > 0$ becomes another parameter.

Here we now have $\Omega = \{x\in\mathbb{R}:\ 1 - x^2 \geq 0\}$ and $\Gamma_R = \emptyset$, meaning that both Theorem~\ref{thm:PDRConv} and Theorem~\ref{thm:PDRsos} are fully applicable. Furthermore, we may anticipate the nearing singularity at $x = \pm 1$ of the optimizing sequence $f^k(x)$ in Theorem~\ref{thm:PDRConv} by introducing a rational ansatz 
\begin{equation}
    f(x) = \frac{\tilde{f}(x)}{1 - x^2}
\end{equation}
for $f$ in the inequality in \eqref{PDR}. The result is the equivalent polynomial inequality on $(x,y,z) \in \Omega\times\mathbb{R}\times\mathbb{R}$ given by
\begin{equation}\label{DoubleWellINeq}
    \bigg[\frac{1}{L^2}z^2 + (U(Lx) - \lambda)y^2\bigg](1 - x^2)^2 + (1 - x^2)\tilde{f}'(x)y^2 + 2x\tilde{f}(x)y^2 + 2(1- x^2)\tilde{f}(x)yz \geq 0,  
\end{equation}
after multiplying through by $(1 - x^2)^2$. We take $\tilde{f}(x)$ to be a degree $\nu$ polynomial in $x$ and optimize over $\lambda \in \mathbb{R}$, leading to the SOS-constrained polynomial optimization problem to bound $\lambda_0$ from below for any $L > 0$.

\begin{figure}[t] 
\center
\includegraphics[width = 0.8\textwidth]{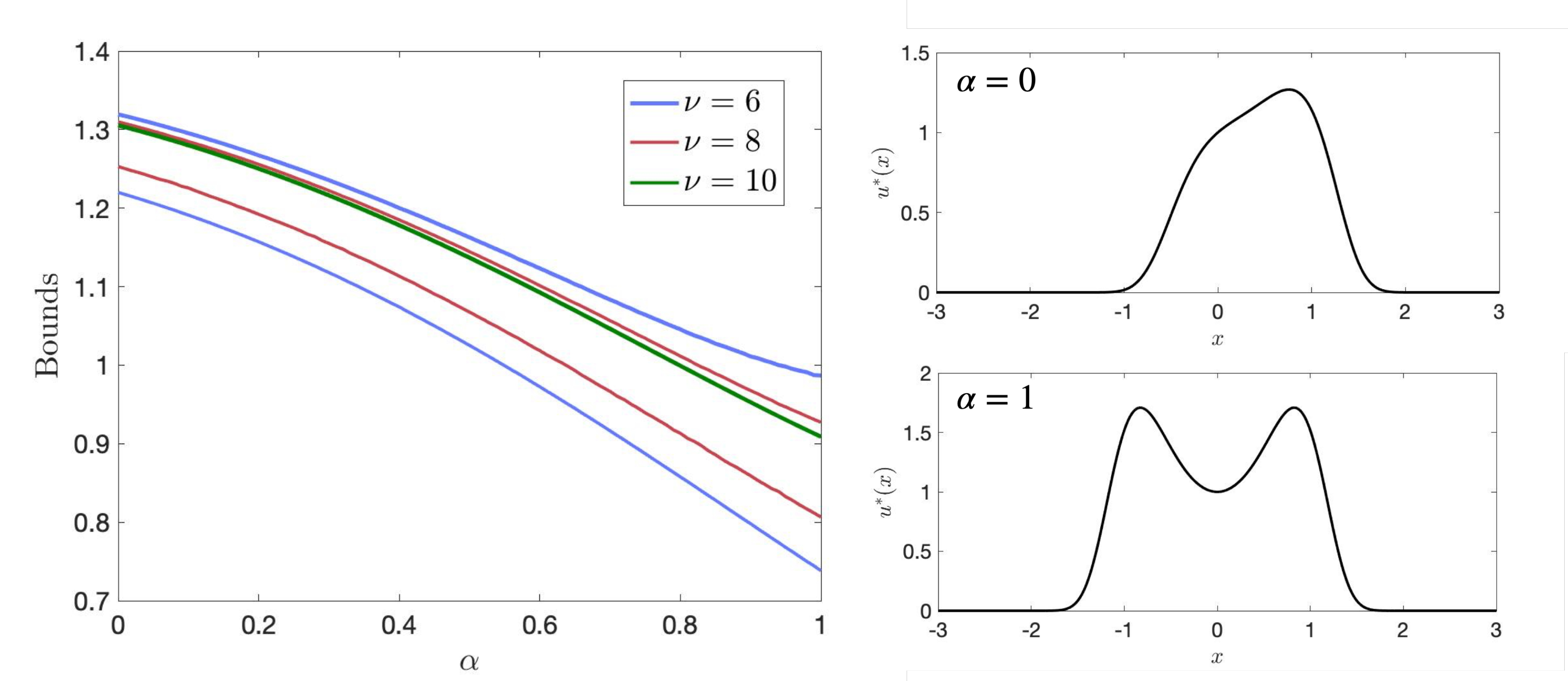}  
\caption{Left: Numerical upper and lower bounds on the principal eigenvalue of the operator \eqref{DoubleWellOperator} with $L = 3$ using degree $\nu = 6,8,10$ polynomials. Right: Approximations of the principal eigenfunction of the operator \eqref{DoubleWellOperator} for $\alpha = 0,1$ using the degree $10$ polynomials that provide the upper bounds at these parameter values.}
\label{fig:Double_Well_Bounds}
\end{figure}

For demonstration we will fix $L = 3$, while noting that larger values of $L$ return similar estimates on the minimal eigenvalue since the principal eigenfunction is localized about the two potential wells at $x = -\alpha,1$. Notice that since $\Omega$ is defined by a single polynomial inequality, we require only on auxiliary SOS polynomial $\rho_1(x,y,z)$ to implement \eqref{PDRsos}. Since the inequality \eqref{DoubleWellINeq} has a maximal degree of 2 in the variables $(y,z)$, the unknown $\rho_1$ can also be taken to be quadratic in these variables as well, significantly decreasing the number of unknown coefficients to be tuned by the numerical procedure \footnote{\color{black} Doing this will not affect the convergence properties of Theorem~\ref{thm:PDRsos}. This is because the proof of Theorem~\ref{thm:PDRsos} in \cite[Appendix B]{chernyavsky2023convex} shows that the weighted SOS representation coming from $\mathcal{Q}_\nu(\Omega)$ need only include polynomials up to the highest degree in $(y,z)$ as $\sigma|z|^2 + (U(x) - \lambda )y^2 + (\nabla\cdot f(x)) y^2 + 2yf(x)\cdot z$.}. Figure~\ref{fig:Double_Well_Bounds} presents numerical upper and lower bounds over $\alpha \in [0,1]$ using the methods of Section~\ref{sec:Bounds} with $\tilde{f}$ of degree $\nu = 6,8,10$. With $\nu = 10$ upper and lower bounds coincide to plotting accuracy, agreeing up to three significant digits, as demonstrated in Table~\ref{tbl:Double_Well}. Figure~\ref{fig:Double_Well_Bounds} further provides approximations of the principal eigenfunctions of $\mathcal{L}$ for $\alpha = 0,1$ obtained from the degree 10 near-optimizers $f(x) = \tilde{f}(x)/(1 - x^2)$. 

\begin{table}[t]
	\centering
	\caption{Upper and lower bounds on the principal eigenvalue of operator \eqref{DoubleWellOperator} with $L = 3$ at different $\alpha\in [0,1]$. Bounds are attained with degree $10$ polynomials and are plotted in Figure~\ref{fig:Double_Well_Bounds}.}
    \begin{tabular}{ r |cccccc}
		\toprule
		$\alpha$ & $0.0$ & $0.2$ & $0.4$ & $0.6$ &  $0.8$ & $1.0$  \\ [0.5ex]
		\midrule
		Upper Bound & 1.3061 & 1.2509 & 1.1790 & 1.0933 & 0.9999 & 0.9095  \\
        Lower Bound & 1.3047 & 1.2495 & 1.1775 & 1.0918 & 0.9984 & 0.9082 \\
		\bottomrule
	\end{tabular}
	\label{tbl:Double_Well}
\end{table}

\subsection{Brownian Motion in a Ball}

In this demonstration we consider the motion of a Brownian particle in a ball in $\mathbb{R}^2$ and $\mathbb{R}^3$. Precisely, the governing equations are 
\begin{equation}
    \mathrm{d}X(t) = \sqrt{2}\mathrm{d}B(t),
\end{equation}
with $X(t) \in \{x: \|x\| < 1\}$ and notably $\nabla V \equiv 0$ here. Such a problem subject to various boundary conditions on $\partial \Omega$ is well-studied in the mathematical literature \cite{cheviakov2010asymptotic,pillay2010asymptotic,delgado2015conditional,iyaniwura2021simulation,lelievre2024spectral,banuelos2023bounds}. This is due to the problem's application in biophysics and cell biology wherein an ion modeled by a Brownian particle attempts to escape a biological cell through membrane channels. The reader is directed to the works \cite{grebenkov2019full,holcman2004escape,schuss2007narrow,grigoriev2002kinetics,benichou2008narrow} and the references therein for a proper overview of the application of the problem to biology. From the mathematical perspective we note that the generator is $\mathcal{L} = \Delta$, the Laplacian operator on $\Omega$, which is already in symmetrized form, i.e. $\mathcal{H} = \mathcal{L}$.

\subsubsection{Fully Absorbing Boundary}\label{subsec:FullyAbsorbing}

We begin with a simple demonstration of the performance of the method to cases where the principal eigenvalue is known. The goal of this brief presentation is to highlight key aspects of the numerical implementation for the more complex demonstration that follow for partially absorbing boundaries. Here we will consider full Dirichlet boundary conditions on the unit ball $\Omega = \{x:\ \|x\| < 1\}$ in $\mathbb{R}^2$ and $\mathbb{R}^3$. As in the previous double well demonstration, we have $\Gamma_A = \partial \Omega$ and $\Gamma_R = \emptyset$, meaning that Theorems~\ref{thm:PDRConv} and \ref{thm:PDRsos} fully apply.   

Start with $\Omega \subset \mathbb{R}^2$ for exposition. The $x_1 \mapsto -x_1$ and $x_2 \mapsto -x_2$ symmetry of the unit ball in $\mathbb{R}^2$ combined with the equivariance of the Laplacian with respect to these symmetries guarantees that the unique principal eigenfunction is such that $u_0(x_1,x_2)$ is even in both of its arguments\footnote{\color{black}One can take this argument further to find that $u_0$ is radially symmetric and incorporating such radial symmetry would re-cast the Laplacian as $\Delta = \partial^2_r + r^{-1}\partial_r$, where $r = \|x\|$. Doing so here would require symmetrizing the eigenproblem to an associated Witten Laplacian or working with an integrand of $r|\partial_r u|^2 - \lambda r u^2$ in \eqref{RQ}. Since all terms in the later formulation are polynomial, the methods herein can be applied to this setting but we refrain from doing so since the results are strong enough with only the assumed reflection symmetries.}. Then, since the optimizing sequence in Theorem~\ref{thm:PDRConv} is approaching $\nabla u_0/u_0$, we consider $f = (f_1,f_2)$ so that $f_1$ is odd in $x_1$ and even in $x_2$, while $f_2$ is even in $x_1$ and odd in $x_2$\footnote{\color{black}Restricting our search to such symmetric polynomials will not affect the convergence results in Theorems~\ref{thm:PDRsos} as one can always find symmetric polynomials that result in the same bound $\lambda_\nu$ in \eqref{PDRsos} for a given degree $\nu$; see \cite[Proposition~6]{oeri2023convex} or \cite[Proposition~4.1]{fantuzzi2020bounding} for methods of proving this.}. Furthermore, we again anticipate the singularity at the boundary of $\Omega$ coming from the Dirichlet boundary conditions by introducing a rational ansatz $f = \tilde{f}/(1 - x_1^2 - x_2^2)$. Thus, the polynomial inequality that defines the problem \eqref{PDR} is given by
\begin{equation}\label{2dDirichlet}
    \begin{split}
    [z_1^2 + z_2^2 - &\lambda y^2](1 - x_1^2 - x_2^2)^2 + \frac{\partial \tilde f_1}{\partial x_1}(x_1,x_2)(1 - x_1^2 - x_2^2)y^2 + 2x_1\tilde f_1(x_1,x_2)y^2 \\
    & + \frac{\partial \tilde f_2}{\partial x_2}(x_1,x_2)(1 - x_1^2 - x_2^2)y^2 + 2x_2\tilde f_1(x_1,x_2)y^2 \\ 
    & + 2(\tilde f_1(x_1,x_2)z_1 + \tilde f_2(x_1,x_2)z_2)y(1 - x_1^2 - x_2^2) \geq 0  
    \end{split}
\end{equation}
where $z_i$ take the place of $\partial u/\partial x_i$ for $i = 1,2$. Again $\Omega = \{1 - x_1^2 - x_2^2 \geq 0\}$ is defined by a single polynomial inequality, necessitating only $\rho_1(x_1,x_2,y,z_1,z_2)$ to implement the SOS problem \eqref{PDRsos}. Moreover, the symmetries of $f = (f_1,f_2)$ endow \eqref{2dDirichlet} with invariances with respect to the actions $(x_1,z_1) \mapsto -(x_1,z_1)$ and $(x_2,z_2) \mapsto -(x_2,z_2)$. Thus, we may take $\rho_1$ to respect these invariances as well, while further restricting it to have no more than quadratic degree in $(y,z_1,z_2)$ since this is also the case for \eqref{2dDirichlet}. In implementing these restrictions we greatly reduce the number of variables/coefficients in the SOS program which can help to mitigate poor numerical conditioning that can arise when solving large SDPs.   

\begin{table}[t]
	\centering
	\caption{Lower bounds on the principal Dirichlet eigenvalue of the Laplacian operator on the unit ball in 2D (left) and 3D (right). Exact values for the principal eigenvalue are 5.7832 in 2D and 9.8696 in 3D.}
	\begin{tabular}{ c |c }
		\toprule
		Degree & Lower Bound  \\ 
		\midrule
        4 & 5.7393 \\
        6 & 5.7811  \\
        8 & 5.7832  \\
        10 & 5.7832 \\
        12 & 5.7832
	\end{tabular}\qquad\qquad 
 \begin{tabular}{ c |c}
		\toprule
		Degree & Lower Bound   \\ 
		\midrule
        4 &  9.7564\\
        6 &  9.8613\\
        8 &  9.8694\\
        10 & 9.8696 \\
        12 & 9.8696
	\end{tabular}
	\label{tbl:DirichletBall}
\end{table}

The implementation for $\Omega$ being unit ball in 3D is nearly identical, but exhibits an additional $x_3 \mapsto -x_3$ symmetry that is reflected in the polynomials. Table~\ref{tbl:DirichletBall} presents lower bounds with increasing polynomial degree of $f$ on the principal eigenvalues that quickly coincide at all presented digits. Moreover, the exact value of the principal eigenvalues can be identified using separation of variables as approximately 5.7832 in 2D and exactly $\pi^2 \approx 9.8696$ in 3D. Upper bounds quickly match the exact value and so are not presented. For example, at degree 4 in 2D the upper bounds match the exact value to three decimal places. Such high accuracy of the upper bounds for low polynomial degree was also observed in Figure~\ref{fig:Double_Well_Bounds} in the previous demonstration.

\subsubsection{Partially Absorbing Boundary}\label{subsec:NarrowEscape}

Let us now consider an escape problem wherein most of the boundary of a ball in 2D is reflecting with only small exits from which the random walker can escape. We begin with a general setting for the problem by considering a sequence of $k \geq 1$ points on the boundary $x_b^{(j)} \in \Omega = \{(x_1,x_2):\ x_1^2 + x_2^2 = 1\}$ and a radius $r \geq 0$ to define 
\begin{equation}\label{EscapeLocations}
    \Gamma_A = \bigcup_{j = 1}^k\ \{x\in\mathbb{R}^2:\ \|x - x_b^{(j)}\| \leq r\ \mathrm{and}\ x \in \partial\Omega\}, 
\end{equation}
while the reflecting portion is what remains of the boundary, i.e. $\Gamma_R = \partial\Omega\setminus\Gamma_A$. Notice that when $r = 0$ the problem reduces entirely to identifying Laplacian eigenvalues on the unit ball with Neumann boundary conditions, meaning that the principal eigenvalue is $\lambda_0 = 0$ with constant eigenfunction. With $r > 0$ the work \cite{lelievre2024spectral} shows that not only is the principal eigenvalue $\lambda_0$ positive, but in the asymptotic regime $0 < r \ll 1$ we have $\lambda_0 \sim -1/\log(r)$. The goal of this demonstration is to show that the methods herein extend these results into the non-asymptotic regime for moderate values of $r$ that is not necessarily amenable to analysis.

\begin{figure}[t] 
\center
\includegraphics[width = 0.8\textwidth]{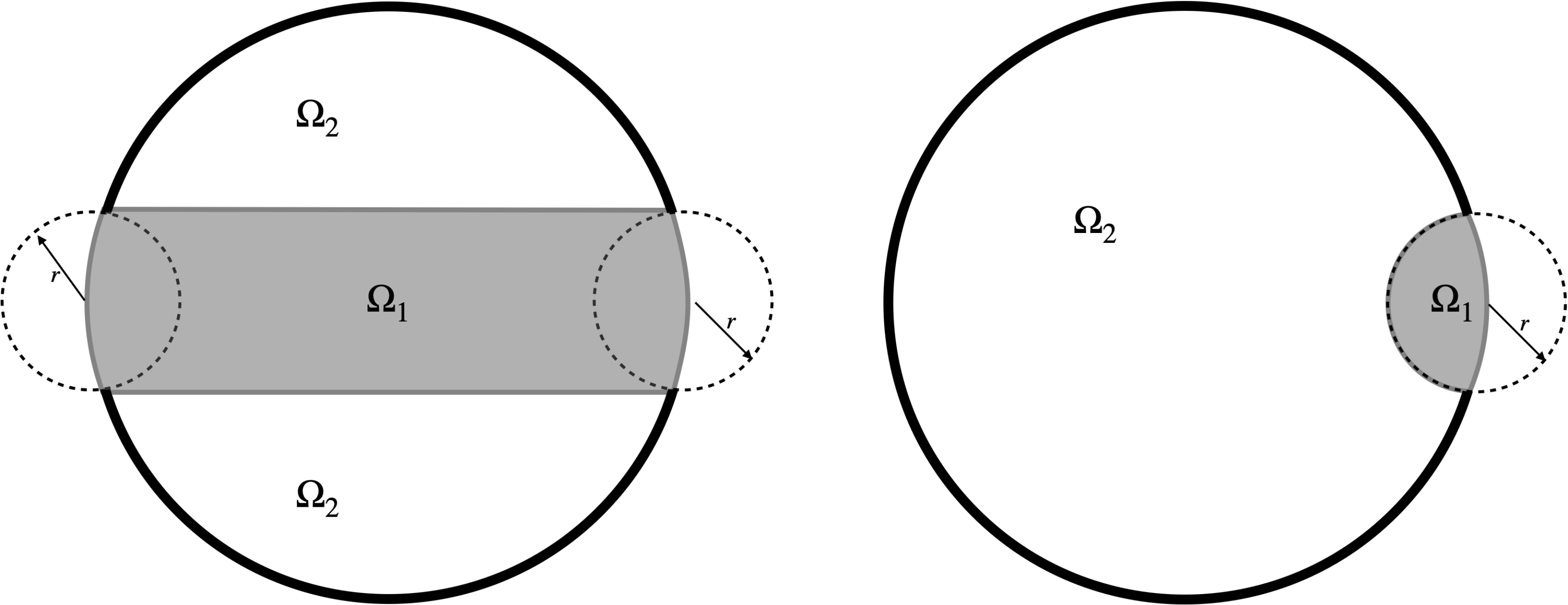}  
\caption{Left: The unit ball in $\mathbb{R}^2$ with exits as in \eqref{EscapeLocations} symmetrically centered at $(\pm 1,0)$ on the boundary. Right: A single exit centered at $(1,0)$ on the boundary.}
\label{fig:Narrow_Escape}
\end{figure}

\paragraph{Two Exit Regions.} Let us consider first the case $k = 2$ with exits along the boundary centered at $(\pm 1,0)$. To implement the bounding method for this problem, we decompose the domain into disjoint subsets $\Omega_1$ and $\Omega_2$ so that their closure intersects $\partial\Omega$ at only $\Gamma_A$ or $\Gamma_R$, but not both. Precisely, let $b = \sup\{x:\ x \in \Gamma_A\} = \frac{r}{2}\sqrt{(4 - r^2))}$ and consider the semialgebraic sets 
\begin{equation}
    \begin{split}
        \Omega_1 &= \{(x_1,x_2):\ 1 - x_1^2 - x_2^2 \geq 0, \ b^2 - x_2^2 \geq 0\} \\
        \Omega_2 &= \{(x_1,x_2):\ 1 - x_1^2 - x_2^2 \geq 0, \ x_2^2 - b^2 \geq 0\}
    \end{split}   
\end{equation}
whose union is all of the unit ball $\Omega$. This decomposition of the full domain $\Omega$ is illustrated in the left panel of Figure~\ref{fig:Narrow_Escape}. Notice that $\partial \Omega_1\cap \partial\Omega = \Gamma_A$ and $\partial \Omega_2\cap \partial\Omega = \Gamma_R$ in this case. We now implement the polynomial optimization framework on $\Omega_1$ and $\Omega_2$ separately, but impose a continuity condition at the interfaces between these sets. To begin, the Dirichlet boundary conditions on $\partial \Omega_1$ again lead to a rational ansatz for the null Lagrangian, $f^{(1)}(x_1,x_2) = \tilde{f}^{(1)}(x_1,x_2)/(1 - x_1^2 - x_2^2)$, and so for $(x,y,z) \in \Omega_1\times\mathbb{R}\times\mathbb{R}^2$ we have 
\begin{equation}\label{Omega1}
    \begin{split}
        [z_1^2 + z_2^2 - &\lambda y^2](1 - x_1^2 - x_2^2)^2 + \frac{\partial \tilde f_1^{(1)}}{\partial x_1}(x_1,x_2)(1 - x_1^2 - x_2^2)y^2 + 2x_1\tilde f_1^{(1)}(x_1,x_2)y^2 \\
        &+ \frac{\partial \tilde f_2^{(1)}}{\partial x_2}(x_1,x_2)(1 - x_1^2 - x_2^2)y^2 + 2x_2\tilde f_1^{(1)}(x_1,x_2)y^2 \\
        &+ 2(\tilde f_1^{(1)}(x_1,x_2)z_1 + \tilde f_2^{(1)}(x_1,x_2)z_2)y(1 - x_1^2 - x_2^2) \geq 0.  
    \end{split}
\end{equation}
The boundary conditions on $\Omega_2$ do not necessitate a rational ansatz and so for $(x,y,z) \in \Omega_2\times\mathbb{R}\times\mathbb{R}^2$ we have
\begin{equation}\label{Omega2}
    \begin{split}
        &[z_1^2 + z_2^2 - \lambda y^2] + \frac{\partial f_1^{(2)}}{\partial x_1}(x_1,x_2)y^2 + 2x_1 f_1^{(2)}(x_1,x_2)y^2 + \frac{\partial f_2^{(2)}}{\partial x_2}(x_1,x_2)y^2 \\
        &\quad + 2x_2 f_1^{(2)}(x_1,x_2)y^2 + 2(f_1^{(2)}(x_1,x_2)z_1 + f_2^{(2)}(x_1,x_2)z_2)y \geq 0.  
    \end{split}
\end{equation}
However, the Neumann boundary conditions on $\partial \Omega_2$ now necessitate the inequality conditions
\begin{equation}\label{NeumannBCs}
    \begin{split}
        2x_1f_1^{(2)}(x_1,x_2) + 2x_2f_2^{(2)}(x_1,x_2) &\geq 0 \\
        -2x_1f_1^{(2)}(x_1,x_2) - 2x_2f_2^{(2)}(x_1,x_2) &\geq 0 \\
    \end{split}
\end{equation}
for all $(x_1,x_2) \in \partial \Omega_2 \cap \partial \Omega = \{1 - x_1^2 - x_2^2 = 0,\ x_2^2 - b^2 \geq 0\}$, where we are enforcing the equality $\vec{n}\cdot f^{(2)} = 0$ through polynomial inequality conditions. Finally, we require a continuity condition from $f^{(1)}$ in $\Omega_1$ to $f^{(2)}$ in $\Omega_2$ at their shared boundaries $x_2 = \pm b$, so for all $-1 + \frac{1}{2}r^2 \leq x_1 \leq 1 - \frac{1}{2}r^2$ we impose
\begin{equation}\label{ContinuityConditions}
    \begin{split}
        \tilde{f}^{(1)}(x_1,b) - f^{(2)}(x_1,b)(1 - x_1^2 - b^2) \geq 0 \\
        -\tilde{f}^{(1)}(x_1,b) + f^{(2)}(x_1,b)(1 - x_1^2 - b^2) \geq 0 \\
        \tilde{f}^{(1)}(x_1,-b) - f^{(2)}(x_1,-b)(1 - x_1^2 - b^2) \geq 0 \\
        -\tilde{f}^{(1)}(x_1,-b) + f^{(2)}(x_1,-b)(1 - x_1^2 - b^2) \geq 0
    \end{split}
\end{equation}
where we recall the rational ansatz for $f^{(1)}$ whose denominator gets multiplied through the inequality conditions to guarantee everything is polynomial when $f^{(1)}$ and $f^{(2)}$ are. Together the inequality conditions \eqref{Omega1}, \eqref{Omega2}, \eqref{NeumannBCs}, and \eqref{ContinuityConditions} imposed on their respective semialgebraic domains and optimized over $f^{(1)},f^{(2)}$ and $\lambda$ provide bounds from below on the principal eigenvalue in this case. As in the previous demonstrations, these inequality conditions are strengthened to SOS conditions to be solved as an SDP. Finally, the $x_1 \mapsto -x_1$ and $x_2 \mapsto -x_2$ symmetries of the domains allow one to impose the same symmetry conditions on $f^{(1)}$ and $f^{(2)}$ as were imposed on the null Lagrangian in \S~\ref{subsec:FullyAbsorbing}. 

\begin{table}[t]
	\centering
	\caption{Lower bounds on the principal eigenvalue of the Laplacian operator with $\Gamma_A$ as in \eqref{EscapeLocations} and $k = 2$ exits symmetrically placed at $(\pm 1,0)$ and $r = 0.5$ (left) and $r = 0.4$ (right). Upper bounds computed using degree 20 polynomials are 1.711 (left) and 1.403 (right). Lower bounds for degrees larger than those presented are not reported as the solver failed.  }
	\begin{tabular}{ c |c }
		\toprule
		Degree & Lower Bound  \\ 
		\midrule
        10 & 1.398 \\ 
        20 & 1.676 \\ 
        30 & 1.692 \\ 
        40 & 1.702 \\ 
        50 & 1.708 \\ 
	\end{tabular}\qquad\qquad  
 \begin{tabular}{ c |c}
		\toprule
		Degree & Lower Bound   \\ 
        \midrule
        10 & 0.881 \\
        20 & 1.372  \\
        30 & 1.393  \\
	\end{tabular}
	\label{tbl:2Holes}
\end{table}

\begin{figure}[t] 
\center
\includegraphics[width = 0.8\textwidth]{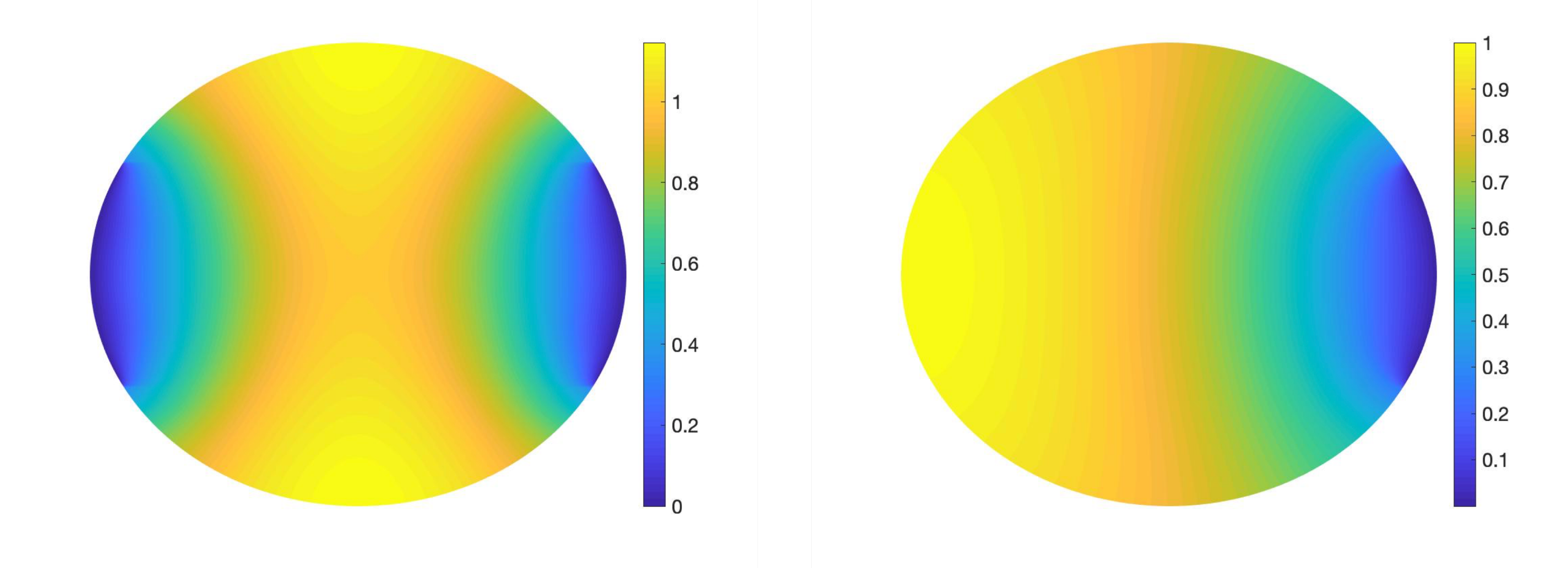}  
\caption{Surface plots of the approximations of the principal eigenfunctions for the Laplacian with partially absorbing boundaries given by \eqref{EscapeLocations} with $(k,r) = (2,0.5)$ (left) and $(k,r) = (1,0.5)$ (right).}
\label{fig:NarrowEscape_Eigenfunctions}
\end{figure}

Table~\ref{tbl:2Holes} presents lower bounds on the principal eigenvalue for this problem with $r = 0.5$ and $r=0.4$ to the trusted precision of three decimal places. Notice that the polynomial degrees on the null Lagrangians $f^{(1)}$ and $f^{(2)}$ are significantly higher than in the previous demonstrations, coming from the complexity of the problem, the polynomial inequality conditions, and the leading eigenfunction. Furthermore, the solver fails at degrees 40 and 50 for $r = 0.4$ and so these lower bounds are not reported in the table. Nonetheless, upper bounds are computed using degree 20 polynomials as 1.711 for $r = 0.5$ and 1.403 for $r = 0.4$, agreeing with the largest lower bounds to at least two significant digits. An approximation of the leading eigenfunction is presented in Figure~\ref{fig:NarrowEscape_Eigenfunctions} for $r = 0.5$, with the $r = 0.4$ case looking similar. 

The failure at high degrees for $r = 0.4$ becomes even more pronounced as $r$ decreases since the eigenfunction becomes singular as $r \to 0^+$. While the singular regime can be amenable to mathematical analysis, it is not suitable for these numerical techniques. For example, decreasing to $r = 0.25$ one finds lower bounds up to degree 30 again, but now the gap between upper and lower bounds is slightly larger, putting the principal eigenvalue in the interval $[0.913,0.938]$. With $r = 0.1$ the region $\Omega_1$ becomes numerically small and appears to have little influence on the lower bounds since even at degree 50 we achieve a lower bound of 0.003. This value is only slightly above the value 0.000 that comes from full Neumann boundary conditions ($\Gamma_A = \emptyset$) and a significant distance from the value of $\lambda_0 \approx 0.76$ reported in \cite{lelievre2024spectral}. 

\paragraph{One Exit Region.} For further demonstration, let us consider the case of only one boundary exit ($k = 1$) located at $(1,0)$. To implement the numerical bounding procedure we again break the domain into two regions, this time given by
\begin{equation}
    \begin{split}
        \Omega_1 &= \{(x_1,x_2):\ 1 - x_1^2 - x_2^2 \geq 0, \ r^2 - x_1^2 - x_2^2 \geq 0\} \\
        \Omega_2 &= \{(x_1,x_2):\ 1 - x_1^2 - x_2^2 \geq 0, \ x_1^2 + x_2^2 - r^2 \geq 0\},
    \end{split}   
\end{equation}
as depicted in the right panel of Figure~\ref{fig:Narrow_Escape}. We again have the inequalities \eqref{Omega1} on the set $\Omega_1$ and \eqref{Omega2} on $\Omega_2$, while \eqref{NeumannBCs} remains but is now enforced on the boundary 
\begin{equation}
    (x_1,x_2)\in\partial\Omega_2 \cap \partial\Omega = \{(x_1,x_2):\ 1 - x_1^2 - x_2^2 = 0, \ x_1^2 + x_2^2 - r^2 \geq 0 \}.    
\end{equation} 
The major difference from the $k = 2$ case above is that the continuity conditions between $f^{(1)}$ and $f^{(2)}$ are enforced on $\partial\Omega_1\cap\partial\Omega_2 = \{(x_1,x_2):\ 1 - x_1^2 - x_2^2 \geq 0,\ r^2 - x_1^2 - x_2^2 = 0\}$ and defined by the inequalities
\begin{equation}\label{ContinuityConditions2}
    \begin{split}
        \tilde{f}^{(1)}(x_1,x_2) - f^{(2)}(x_1,x_2)(1 - x_1^2 - x_2^2) \geq 0 \\
        -\tilde{f}^{(1)}(x_1,x_2) + f^{(2)}(x_1,x_2)(1 - x_1^2 - x_2^2) \geq 0.
    \end{split}
\end{equation}
Finally, since the domains $\Omega$, $\Omega_1$, and $\Omega_2$ no longer exhibit an $x_1 \mapsto -x_1$ invariance, there is no symmetry in $x_1$ that can be imposed on the functions $f^{(1)}$ and $f^{(2)}$. The symmetry in $x_2$ remains. 

\begin{table}[t]
	\centering
	\caption{Lower bounds on the principal eigenvalue of the Laplacian operator with $\Gamma_A$ as in \eqref{EscapeLocations} and $k = 1$ exit placed at $(\pm 1,0)$ and $r = 0.5$ (left) and $r = 0.4$ (right). Upper bounds computed using degree 20 polynomials are 0.563 (left) and 0.628 (right). Lower bounds for degrees larger than those presented are numerically inaccurate as reported by the solver. } 
	\begin{tabular}{ c |c }
		\toprule
		Degree & Lower Bound  \\ 
		\midrule
        10 & 0.340 \\ 
        20 & 0.544 \\ 
        30 & 0.552 \\ 
        40 & 0.556 \\ 
        50 & 0.556 
	\end{tabular}\qquad\qquad 
 \begin{tabular}{ c |c}
		\toprule
		Degree & Lower Bound   \\ 
		\midrule
        10 & 0.092\\ 
        20 & 0.461\\ 
        30 & 0.484 \\ 
        40 & 0.489\\
	\end{tabular}
	\label{tbl:1Hole}
\end{table}  

Numerical lower bounds are presented in Table~\ref{tbl:1Hole} for $r = 0.4,0.5$ and increasing polynomial degrees. Furthermore, an approximation of the leading eigenfunction in the case $r = 0.5$ is presented in Figure~\ref{fig:NarrowEscape_Eigenfunctions}. Again, results are promising, but break down as $r$ becomes small and enters the singular parameter regime. Furthermore, while the degree 20 upper bound of $0.563$ for $r=0.5$ matches the greatest lower bound reported in Table~\ref{tbl:1Hole} to two decimal places, the degree 20 upper bound of $0.627$ for $r = 0.4$ remains distant from the greatest lower bound at $0.489$ reported in Table~\ref{tbl:1Hole}. Increasing the polynomial degree to 30 lowers the upper bound to $0.580$, but we still see a significant break from all previous demonstrations where relatively low polynomial degrees provided tight upper bounds. While we do not offer a hypothesis for the cause of this, we conclude by noting that we have put the principal eigenvalue for $(k,r) = (1,0.4)$ into the interval $[0.489,0.580]$ which is still quite small.

\section{Discussion}\label{sec:Discussion}

In this work we have described and applied a numerical method for bounding principal eigenvalues of generators of stochastic processes from both above and below. Our focus here was on Brownian dynamics, but the method is general enough to be applied to a wider class of stochastic differential equations. The applications in Section~\ref{sec:Demonstrations} demonstrated the applicability on the model problem of escaping from a double well potential landscape, as well as on the physically-relevant problem of a Brownian particle escaping a ball. Moreover, we saw that these {\color{black} numerical} upper and lower bounds are complemented with an approximation of the QSD, i.e. the leading eigenfunction. While the method as it is presented in this manuscript is only applicable to polynomial potential functions, one can apply various tricks to re-cast rational or trigonometric nonlinearities into polynomial form (see for example \cite{bramburger2024synthesizing}). Furthermore, one may additionally consider state-dependent $\sigma$ with no hindrance to the application of the method so long as the existence of an isolated leading eigenvalue of the generator is known a priori. The result is a robust {\color{black} numerical} method that is firmly rooted in theory which can complement pencil-and-paper analysis to provide a better understanding of the probabilistic behaviour of stochastic differential equations.   

A potential avenue for further application of this method could be to large deviations of Markov processes. Indeed, the G\"artner--Ellis theorem relates large deviations of stochastic processes to the leading eigenfunction of a {\em tilted generator} \cite{touchette2018introduction}. Thus, at its core, one can apply the methods presented in this paper to bound these leading eigenvalues from above and below. While recent work has employed computer-assisted proofs to rigorously bound these leading eigenvalues to detect random bifurcations in stochastic systems \cite{blessing2024detecting}, the SOS methods herein provide an alternative, albeit non-rigorous, method to quickly and efficiently bound these eigenvalues to good accuracy. Furthermore, the bounds and approximation of the leading eigenfunction from these SOS programs could function as a pre-conditioner for rigorous numerical routines by providing accurate initial guesses to start a computer-assisted proof.

{\color{black}Finally, it should be emphasized that finding SOS representations using SDPs can be made rigorous using rational arithmetic \cite{magron2018exact,henrion2016exact}, interval arithmetic \cite{jansson2008rigorous}, or high precision arithmetic \cite{habibi2024loraine,leijenhorst2024solving}. This means that one could provide validated upper and lower bounds on $\lambda_0$ with such numerical methods, thereby further bridging the gap between theory and implementation.} What is most important for any such SOS program is that it is implemented in such a way to keep the problem numerically well-conditioned. This means building in any and all symmetries to reduce the number of polynomial coefficients, as well as incorporating analytical information to reduce the computational complexity of the problem, as was done with the rational ansatz in our demonstrations. Without such information it is often difficult to achieve sharp bounds since numerical issues appear at high polynomial degrees due to the size and complexity of the SDP to be solved. ackThese problems become especially apparent for even moderate dimensions of $x \in \Omega\subset \mathbb{R}^d$, thus making this method in its current form best-suited for low-dimensional problems.

\section*{Acknowledgment}

This work was partially supported by an NSERC Discovery Grant and the Fondes de Recherche du Qu\'ebec – Nature et Technologies (FRQNT). The author gratefully acknowledges helpful comments and discussions with Maximilian Engel, Giovanni Fantuzzi, and David Goluskin. A special acknowledgment is made to Tony Leli\`evre for providing discussions and helpful comments on drafts of this manuscript. Finally, the author would like to thank the anonymous referees for helping to improve the paper with their constructive comments.

\section*{Data Availability Statement}

All code to reproduce the numerical results reported in this paper are publicly available for download and use at: \href{https://github.com/jbramburger/SDE-Escape-Rates}{https://github.com/jbramburger/SDE-Escape-Rates}.

\bibliographystyle{abbrv}
\bibliography{references.bib}

\end{document}